\documentclass[12pt, reqno]{amsart}
\usepackage{amsmath, amsthm, amscd, amsfonts, amssymb, graphicx, xcolor}
\usepackage[bookmarksnumbered, colorlinks, plainpages]{hyperref}
\usepackage{amsmath}
\usepackage{cases}
\usepackage{amssymb}
\usepackage{mathrsfs}
\usepackage{euscript}
\usepackage{dsfont}
\usepackage{graphicx}
\usepackage{float}
\usepackage{graphicx}
\usepackage{amssymb}

\allowdisplaybreaks[4]
\newtheorem{theorem}{Theorem}[section]

\newtheorem{lemma}{Lemma}[section]

\theoremstyle{definition}

\theoremstyle{remark}
\newtheorem*{remark}{Remark}

\numberwithin{equation}{section}

\begin{document}
	\setcounter{page}{1}
	

	
	\title[]{Note on fractional sum of the divisor function}
	
	\author[]{Zhiyuan Yang}
	
	\address{School of Mathematics,  Shandong University, Jinan 250100, Shandong, China}
	\email{\tt zhiyuan.yang@mail.sdu.edu.cn}

	
	

	
	
	\begin{abstract}
		Let $d(n)$ be the divisor function and denote by $[t]$ the integral part of the real number $t$. We prove that
		\begin{align*}
			\sum_{n\leq x^{1/c}}d\left(\left[\frac{x}{n^c}\right]\right)=d_cx^{1/c}+O_{\varepsilon,c}(x^{\theta_c+\varepsilon}),
		\end{align*}
		where $d_c$ is a suitable constant, $\theta_c<2/(3c+2)$ for $0<c<2/11$.
		This result constitutes an improvement upon that of Y. Feng (2024).
		\newline
		\noindent 
	\end{abstract} \maketitle
	
	
	\section{Introduction}
	In 2019, Bordellès, Dai, Heyman, Pan and Shparlinski [2] studied a general form
	\begin{align*}
		S_f(x):=\sum_{n\leq x}f\left(\left[\frac{x}{n}\right]\right),
	\end{align*}
	under some reasonable hypotheses on $f$. They proved that if there exists $\alpha\in(0,2)$ such that
	\begin{align*}
		\sum_{n\leq x}|f(n)|^2\ll x^\alpha,
	\end{align*}
	then
	\begin{align*}
		S_f(x)=x\sum_{n\geq 1}\frac{f(n)}{n(n+1)}+O\left(x^{(\alpha+1)/3}(\log x)^{(\alpha+1)/3+o(1)}\right).
	\end{align*}
	In [12], Wu improved the error term by establishing the following better results\begin{align*}
		S_f(x)=x\sum_{n\geq 1}\frac{f(n)}{n(n+1)}+O\left(x^{(\theta+1)/2}\right),
	\end{align*}
	if there is a constant $\theta\in[0,1)$
	such that $|f(n)|\ll n^\theta$ for $n\geq 1$.
	
	In particular, for the von Mangoldt function $\Lambda(n)$, this is a problem that interests many people. In [10], Ma and Wu firstly break the $1/2$-barrier by proving a better estimate \begin{align*}
		S_{\Lambda}(x)=x\sum_{n\geq 1}\frac{\Lambda(n)}{n(n+1)}+O\left(x^{35/71+\varepsilon}\right).
	\end{align*}
	Latter, the exponent 35/71 was improved to 97/203, 9/19, 92/195, 22/47 by Bordellès
	[1], Liu-Wu-Yang [7], Zhang [14] and Lü-Xu [8] respectively. The best known exponent is $7/15$ due to Chen-Hu-Wu [3].
	
	In this paper, we are interesting in the divisor function $d$. In 2021, Ma and Sun [9] proved that
	\begin{align*}
		S_{d}(x)=x\sum_{n\geq 1}\frac{d(n)}{n(n+1)}+O\left(x^{11/23+\varepsilon}\right).
	\end{align*}
	The exponent 11/23 also was improved to 19/40, 5/11 by Bordellès [1] and Stucky [11] respectively. On the other hand, Feng [4] considered the asymptotic formula of the sum 
	\begin{align*}
		S_{d,c}(x)=\sum_{n\leq x^{1/c}}d\left(\left[\frac{x}{n^c}\right]\right),
	\end{align*}
	where $c>0$. He showed
	\begin{align*}
		S_{d,c}(x)=x^{1/c}\sum_{k\geq 1}d(k)\left(\frac{1}{k^{1/c}}-\frac{1}{(k+1)^{1/c}}\right)+O\left(x^{\theta_c+\varepsilon}\right),
	\end{align*}
	where 
		\begin{equation}
			\begin{aligned}
				\theta_c=&\left\{ \begin{aligned}
					&\frac{2}{3c+2}, \quad&&\text{if}\ 0<c<2/11;\\&\frac{11}{11c+12},\quad&& \text{if}\ c\geq 2/11.
				\end{aligned}\right.\nonumber
			\end{aligned}
		\end{equation}
	Later, motivated by the key ideas used in [11], for $c>2/9$, Wu [13] showed a better error term
	\begin{equation}
		\begin{aligned}
			\theta_c=&\left\{ \begin{aligned}
				&\frac{2c+2}{2c^2+5c+2}, \quad&&\text{if}\ 0<c<2/3;\\&\frac{5}{5c+6},\quad&& \text{if}\ c\geq 2/3.
			\end{aligned}\right.\nonumber
		\end{aligned}
	\end{equation}
	Recently, Feng and Wu [5] obtained a better exponent
	\begin{equation}
		\begin{aligned}
			\theta_c=&\left\{ \begin{aligned}
				&\frac{2}{3c+2}, \quad&&\text{if}\ 0<c<2/5;\\&\frac{5}{5c+6},\quad&& \text{if}\ c\geq 2/5.
			\end{aligned}\right.\nonumber
		\end{aligned}
	\end{equation}
	In this paper, for $0<c<2/13$, we obtain a better error term. Our result are as following.
	\begin{theorem}
		We have
		\begin{align*}
			S_{d,c}=x^{1/c}\sum_{k\geq 1}d(k)\left(\frac{1}{k^{1/c}}-\frac{1}{(k+1)^{1/c}}\right)+O(x^{\theta_c+\varepsilon}),
		\end{align*}
		as $x\rightarrow \infty$, where
		\begin{align*}
				\theta_c=&\min\left(\frac{2}{(3+i)c+2},\frac{1+\kappa_{i+1}+\lambda_{i+1}}{(1+\kappa_{i+1}+\lambda_{i+1})c+2}\right),
		\end{align*}
		for $ 2\kappa_{i+1}/(1+\kappa_{i+1}+\lambda_{i+1})\leq c<2\kappa_i/(1+\kappa_i+\lambda_i)$, $i\geq 0$,
	the pair $(\kappa_i,\lambda_i)$ satisfies
	\begin{align*}
		(\kappa_0,\lambda_0)=\left(\frac{1}{6},\frac{2}{3}\right),\quad (\kappa_{i+1},\lambda_{i+1})=\left(\frac{\kappa_i}{2\kappa_i+2},\frac{\kappa_i+\lambda_i+1}{2\kappa_i+2}\right), \quad i\geq 0.
	\end{align*}
	\end{theorem}
	\begin{remark}
	Note that for $0<c<2/13$, one has $\theta_c<2/(3c+2)$. In fact, this range can be extended at least up to $2/11$.
	\end{remark}

	\section{Lemmas}
	
	Define $e(t):=e^{2\pi it}$.
	The first lemma is following from [13, (3.1)--(3.5)].
	\begin{lemma}
		For $0<c\leq 1$, we have
		\begin{align*}
			S_{d,c}=&x^{1/c}\sum_{k\geq 1}d(k)\left(\frac{1}{k^{1/c}}-\frac{1}{(k+1)^{1/c}}\right)\\&+O\Bigg(x^\varepsilon\left( N+\sum_{\delta=0,1}\max_{N<D\leq x/N^c}\min_{H\geq 1}\mathfrak{S}_{\delta,c}(x,D)\right)\Bigg),
		\end{align*}
		where $1\leq N\leq x^{1/(1+c)}$,
		\begin{align*}
			\mathfrak{S}_{\delta,c}(x,D)=\frac{D}{H}+ \sum_{0<h\leq H}\frac{1}{h}\Bigg|\sum_{D<k\leq 2D}d(k)e\left(\frac{hx^{1/c}}{(k+\delta)^{1/c}}\right)\Bigg|.
		\end{align*}
	\end{lemma}
	
	In this paper, our improvement come from the following so-called exponent pairs. For more details, we can see [6, \S 3. The method of exponent pairs].
	\begin{lemma}
		If $(\kappa, \lambda)$ is an exponent pair, then
		\begin{align*}
			A(\kappa,\lambda)=\left(\frac{\kappa}{2\kappa+2},\frac{\kappa+\lambda+1}{2\kappa+2}\right)
		\end{align*}
		and
		\begin{align*}
			B(\kappa,\lambda)=(\lambda-1/2,\kappa+1/2)
		\end{align*}
		are two exponent pairs. In particular, $(1/6,2/3)$ and $A(1/6,2/3)=(1/14,11/14)$ are two exponent pairs.
	\end{lemma}
	The third lemma is from [6, (3.3.4)].
	\begin{lemma}
		For any exponent pair $(\kappa,\lambda)$, for every $s>0$ and $y>0$, we have
		\begin{align*}
			\sum_{n\in I}e(f(n))\ll (yN^{-s})^\kappa N^\lambda+y^{-1}N^s,
		\end{align*}
		where $I=[a,b]\subset [N,2N]$, the function $f$ satisfies the following. There is some $P=P(\kappa,\lambda,s)$ and some $\varepsilon=\varepsilon(\kappa,\lambda,s)<1/2$ such that
		\begin{itemize}
			\item $f$ is defined and has $P$ continuous derivatives on interval $[a,b]$,
			\item if $0\leq p\leq P-1$ and $a\leq x\leq b$ then
			\begin{align*}
				\Bigg|f^{(p+1)}(x)-(-1)^p(s)_pyx^{-s-p}\Bigg|<\varepsilon (s)_pyx^{-s-p},
			\end{align*}
		\end{itemize}
		where $(s)_0=1$ and $(s)_p=s(s+1)\cdots(s+p-1)$ for $p\geq 1$.
	\end{lemma}

	\section{Proof of Theorem 1.1}
	Note that $d=1*1$.
	By Lemma 2.1, we only need to estimate
	\begin{align*}
		\mathfrak{S}_{\delta,c}(x,D)-\frac{D}{H}&=\sum_{0<h\leq H}\frac{1}{h}\Bigg|\sum_{D<k\leq 2D}d(k)e\left(\frac{hx^{1/c}}{(k+\delta)^{1/c}}\right)\Bigg|\\&=\sum_{0<h\leq H}\frac{1}{h}\Bigg|\sum_{D<mn\leq 2D}e\left(\frac{hx^{1/c}}{(mn+\delta)^{1/c}}\right)\Bigg|\\&\ll\sum_{0<h\leq H}\frac{1}{h}\sum_{m<D^{1/2}}\Bigg|\sum_{D/m<n\leq 2D/m}e\left(\frac{hx^{1/c}}{(mn+\delta)^{1/c}}\right)\Bigg|.
	\end{align*}
	Then we use Lemma 2.3 with
	\begin{align*}
		f(z)=\frac{hx^{1/c}}{(mz+\delta)^{1/c}},\quad s=\frac{1}{c}+1,\quad y=\frac{hx^{1/c}}{cm^{1/c}}
	\end{align*}
	to get
	\begin{align*}
		\sum_{D/m<n\leq 2D/m}e\left(\frac{hx^{1/c}}{(mn+\delta)^{1/c}}\right)&\ll \left(\frac{hx^{1/c}}{m^{1/c}}\frac{1}{(D/m)^{1/c+1}}\right)^{\kappa}\left(\frac{D}{m}\right)^\lambda+\frac{m^{1/c}(D/m)^{1+1/c}}{hx^{1/c}}\\&\ll x^{\kappa/c}D^{-\kappa(1+1/c)+\lambda}m^{\kappa-\lambda}h^\kappa+x^{-1/c}D^{1+1/c}m^{-1}h^{-1}.
	\end{align*}
	Summing over all $h$ and $m$, we get
	\begin{align*}
		\mathfrak{S}_{\delta,c}(x,D)&\ll DH^{-1}+x^\varepsilon\sum_{m<D^{1/2}}\left(x^{\kappa/c}D^{-\kappa(1+1/c)+\lambda}m^{\kappa-\lambda}H^\kappa+x^{-1/c}D^{1+1/c}m^{-1}\right)\\&\ll DH^{-1}+x^{\kappa/c+\varepsilon}D^{-\kappa(1/2+1/c)+\lambda/2+1/2}H^\kappa+x^{-1/c+\varepsilon}D^{1+1/c}.\tag{3.1}
	\end{align*}
	Then we have
	\begin{align*}
		S_{d,c}=x^{1/c}\sum_{k\geq 1}d(k)\left(\frac{1}{k^{1/c}}-\frac{1}{(k+1)^{1/c}}\right)+O\left(x^{2\varepsilon} \left(N+\max_{N<D<x/N^c}R_c(x,D)\right)\right).\tag{3.2}
	\end{align*}
	where
	\begin{align*}
		R_c(x,D)&\ll \min_{H\geq 1}\Big(DH^{-1}+x^{\kappa/c}D^{-\kappa(1/2+1/c)+\lambda/2+1/2}H^\kappa+x^{-1/c}D^{1+1/c}\Big)\\&\ll \min_{H\geq 1}\Big(DH^{-1}+x^{\kappa/c}D^{-\kappa(1/2+1/c)+\lambda/2+1/2}H^\kappa\Big)+xN^{-1-c}.
	\end{align*}
	\subsection{Case 1}
	Let
	\begin{align*}
		&H:=x^{-(\kappa/c)/(\kappa+1)}D^{(\kappa(1/2+1/c)-\lambda/2+1/2)/(\kappa+1)},\\&N:=x^{(\kappa/c)/(\kappa(1/2+1/c)-\lambda/2+1/2)},
	\end{align*}
	thus $H\geq 1$.
	Suppose
	\begin{align*}
		N>xN^{-1-c},\quad N\leq x^{1/(1+c)},\quad -\kappa(-1/2+1/c)+\lambda/2+1/2<0.\tag{3.4}
	\end{align*}
	Then we have
	\begin{align*}
		N+\max_{N<D<x/N^c}R_c(x,D)&\ll N+\max_{N<D<x/N^c}x^{(\kappa/c)/(\kappa+1)}D^{(-\kappa(-1/2+1/c)+\lambda/2+1/2)/(\kappa+1)}\\&\ll N+x^{(\kappa/c)/(\kappa+1)}N^{(-\kappa(-1/2+1/c)+\lambda/2+1/2)/(\kappa+1)}\\&\ll N,
	\end{align*}
	since
	\begin{align*}
		x^{(\kappa/c)/(\kappa+1)}N^{(-\kappa(-1/2+1/c)+\lambda/2+1/2)/(\kappa+1)}\leq N\Leftrightarrow N\geq x^{(\kappa/c)/(\kappa(1/2+1/c)-\lambda/2+1/2)}.
	\end{align*}
	We can rewrite $N$ as
	\begin{align*}
		\frac{\log N}{\log x}&=\frac{\kappa/c}{\kappa(1/2+1/c)-\lambda/2+1/2}\\&=\frac{1/c}{1/2+1/c+(-\lambda+1)/(2\kappa)}.
	\end{align*}
	If $(\kappa,\lambda)=(1/6,2/3)$ with $(-\lambda+1)/(2\kappa)=1$, we have
	\begin{align*}
		\frac{\log N}{\log x}=\frac{1/c}{1/2+1/c+1}=\frac{2}{3c+2},
	\end{align*}
	which is Y. Feng's result [4]. 
	Define $(\kappa_1,\lambda_1):=(1/14,11/14)$ with $(-\lambda_1+1)/(2\kappa_1)=3/2$, and
	\begin{align*}
		A(\kappa_i,\lambda_i)=\left(\frac{\kappa_i}{2\kappa_i+2},\frac{\kappa_i+\lambda_i+1}{2\kappa_i+2}\right):=(\kappa_{i+1},\lambda_{i+1}).
	\end{align*}
	Then we have
	\begin{align*}
		\frac{-\lambda_{i+1}+1}{2\kappa_{i+1}}=-\frac{-(\kappa_i+\lambda_i+1)/(2\kappa_i+2)+1}{2\kappa_i/(2\kappa_i+2)}=\frac{-\lambda_i+1}{2\kappa_i}+\frac{1}{2}.\tag{3.5}
	\end{align*}
	We obtain the result in Theorem 1.1.
	
	Now we need to verify that (3.4) holds.
	The third condition $-\kappa(-1/2+1/c)+\lambda/2+1/2<0$ is equivalent to 
	\begin{align*}
		c<\frac{2\kappa}{1+\kappa+\lambda},
	\end{align*}
	which is implied by the Theorem Hypothesis.
	The first condition $N>xN^{-1-c}$ is equivalent to
	\begin{align*}
		\frac{1/c}{1/2+1/c+(-\lambda+1)/(2\kappa)}>\frac{1}{2+c},
	\end{align*}
	which is implied by
	\begin{align*}
		c<\frac{2\kappa}{1-\kappa-\lambda}.
	\end{align*}
	We also have $N<x^{1/(1+c)}$ since
	\begin{align*}
		\frac{1/c}{1/2+1/c+(-\lambda+1)/(2\kappa)}<\frac{1}{1+c}\Leftrightarrow \frac{-\lambda+1}{2\kappa}>\frac{1}{2}.
	\end{align*}
	We complete the proof of this case.
	\subsection{Case 2} Consider
	\begin{align*}
		-\kappa(-1/2+1/c)+\lambda/2+1/2\geq 0\Leftrightarrow c\geq \frac{2\kappa}{1+\kappa+\lambda}.
	\end{align*}Let\begin{align*}
	&H:=x^{-(\kappa/c)/(\kappa+1)}D^{(\kappa(1/2+1/c)-\lambda/2+1/2)/(\kappa+1)},\\&
	\frac{\log N}{\log x}:=\frac{(1+\lambda+\kappa)/(2c)}{(1+\kappa+\lambda)/2+1/c}.
	\end{align*}
	Then we have $H\ge 1$ and $x^{1/(2+c)}<N<x^{1/(1+c)}$.
	Similar to Case 1, we have
	\begin{align*}
		&N+\max_{N<D<x/N^c}R_c(x,D)\\&\ll N+\max_{N<D<x/N^c}x^{(\kappa/c)/(\kappa+1)}D^{(-\kappa(-1/2+1/c)+\lambda/2+1/2)/(\kappa+1)}\\ &\ll N+x^{(\kappa/c)/(\kappa+1)}(x/N^c)^{(-\kappa(-1/2+1/c)+\lambda/2+1/2)/(\kappa+1)}\\&\ll N.
	\end{align*}
	We also have (we compare $N$ here with that in Case 1 and subtract $1/2$ in the denominator.)
	\begin{align*}
		\frac{(1+\kappa+\lambda)/(2c)}{(1+\kappa+\lambda)/2+1/c}<\frac{1/c}{1/c+(-\lambda+1)/(2\kappa)}\Leftrightarrow c<\frac{2\kappa(1-\kappa-\lambda)}{(1+\kappa+\lambda)(1-2\kappa-\lambda)}.
	\end{align*}
	Then for $(\kappa,\lambda)=(\kappa_{i+1},\lambda_{i+1})$, we have
	\begin{align*}
		\frac{2\kappa_{i+1}(1-\kappa_{i+1}-\lambda_{i+1})}{(1+\kappa_{i+1}+\lambda_{i+1})(1-2\kappa_{i+1}-\lambda_{i+1})}=\frac{2+i}{1+i}\frac{2\kappa_{i+1}}{1+\kappa_{i+1}+\lambda_{i+1}}<\frac{2\kappa_{i}}{1+\kappa_{i}+\lambda_{i}}.
	\end{align*}

	We have completed all the proofs.
	\begin{remark}
		In our proof, we mainly employ the $A$-transformation for exponent pairs. Adding the $B$-transformation may produce better results in some parameter ranges. For example, for $(\kappa,\lambda):=(1/9,13/18)=AB(1/14,11/14)$, corresponding to Case 1 and Case 2, we can choose
		\begin{align*}
			\theta_c:=\frac{4}{7c+4}<\frac{13}{13c+14},
		\end{align*}
		for $4/39<c<4/33$, and
		\begin{align*}
			\theta_c:=\frac{11}{11c+12}<\min\left(\frac{2}{3c+2},\frac{13}{13c+14}\right),
		\end{align*}
		for $4/33<c<2/11$. However, this would involve substantial computations and yield rather complicated results.
	\end{remark}
	
	\section*{Acknowledgements}
	The author thanks Liuying Wu for his explanation regarding the revisions [13]. The author also thanks Wei Zhang for pointing out the omitted reference [5].

\end{document}